\documentclass[12pt]{article}
\usepackage{graphicx}
\usepackage{enumerate}
\usepackage[ruled]{algorithm2e}
\usepackage{amsmath}

\input amssym.def 
\input amssym.tex 

\begin{document}

\title{A Lower Bound for $R(5,6)$}

\author{
Geoffrey Exoo \\
\small Department of Mathematics and Computer Science \\[-0.8ex]
\small Indiana State University \\[-0.8ex]
\small Terre Haute, IN 47809 \\
\small\tt ge@cs.indstate.edu \\
}

\maketitle

\begin{center}
\end{center}

\begin{abstract}
The known lower bound for the the classical Ramsey number $R(5,6)$ is improved from $58$ to $59$.
The method used to construct the graph is a simple variant of computational methods that have been previously
used to construct Ramsey graphs.  The new method uses the concurrent programming features
of the {\em Go} programming language.
\end{abstract}

\section{Introduction}

The classical ($2$-color) Ramsey number $R(s,t)$ is the smallest integer $n$
such that in any $2$-coloring of the edges of the complete graph $K_n$ there
is a monochromatic copy of $K_s$ in color $1$ or of $K_t$ in color $2$.
The reader is referred to Radziszowski's survey on
Small Ramsey Numbers \cite{survey} for basic terminology and
the current status of a host of problems related to Ramsey numbers.

In this paper, a simple modification of a well known
randomized search procedure is used to find a $(5,6)$-coloring of $K_{58}$, thereby
improving the lower bound for $R(5,6)$ to $59$.

\section{The Old Method}

When dealing with diagonal Ramsey numbers, for example $R(5,5)$,
the method shown below can be effective.
It uses two auxiliary procedures, which we describe first.
The first of these prcedures, $count$, counts the total number
of forbidden subgraphs in a coloring.
The second procedure, $ecount$, counts the number of
forbidden subgraphs of a given color and containing a given edge.

The main loop in the algorithm proceeds by shuffling a list of all edges.
It iterates through the list, computing
the $ecount$ for each edge in both colors.  The edge is then assigned the color
that gives the smaller number of forbidden subgraphs.  If case of ties, the edge
is colored randomly.

When dealing with off diagonal Ramsey problems, in particular $R(5,6)$, this method
is less effective.
If we simply count the number of
forbidden subgraphs for each edge, in each color, compare, and then
decide ties randomly (with a probability of $1/2$ for each color), the
algorithm tends to create colorings that are local
minima and have
a large number of $K_5$'s in the first color and no
$K_6$'s in the second color.

Two simple adjustments improve things considerably.
First we assign a larger weight
to monochromatic $K_5$'s than to monochromatic $K_6$'s, the notion being that
monochromatic $K_5$'s are harder to avoid than monochromatic $K_6$'s and
so we must try harder to avoid them.
A weight of about $5/4$ seems to work well.  More precisely, we
assign a penalty of $5$ to a monochromatic $K_5$ in color $1$ and
a penalty of $4$ to a monochromatic $K_6$ in color $2$.

The second adjustment is to change the probabilities used when assigning a random color
to an edge.  In the case of $R(5,6)$, we want to choose color $1$
less often that color $2$.  A probability of approximately $0.47$ seems to
work best.

The main loop of the old algorithm is summarized in the display below.

\vspace{5mm}

\begin{algorithm}[H]
\KwData{$A$ is the $n \times n$ $2$-color adjacency matrix}
\KwData{$L$ is a list of all edges of $K_n$}
\KwData{$weight$ is a list of two elements $[5, 4]$ indexed by color}
\Begin{
  Randomly color $A$ using colors $1$ and $2$ with $P(color\,1) = 0.47$\;
  \Repeat{$ count(A)=0 $}  {
    Shuffle $L$\;
    \For{$e=uv$ in $L$}{
      \For{$c$ in 1,2}{
        Color edge $e$ with color $c$\;
        Compute $escore[c] = weight[c] \cdot ecount(e,c) $\;
      }
      Color edge $e$ with the color $c$ having the minimum $escore$ \;
      In case of equal $escore$s choose color $1$ with $P(color\,1) = 0.47$\;
    }
  }
}
\caption{The old search procedure for finding $(5,6)$-colorings of $K_{n}$}
\end{algorithm}

\section{The New Method}

Recently this author was involved with a project (unrelated to Ramsey numbers)
that used the {\em Go} programming language,
developed at Google.  The project also gave the author unfettered access to a computer with
$256$ cores.
Go has some very nice features that facilitate concurrent programming.
Of course, there is nothing one can do in Go that cannot also be done in C, but
when dealing with concurrency Go is far more convenient.
In addition, the Go compiler generates relatively fast executables.
When doing Ramsey searches, most of the CPU time is spent counting
monochromatic forbidden subgraphs ($K_5$'s and $K_6$'s in our case).
The code for doing these counts in Go is
essentially the same as the code written in C, and is only (approximately) 20 percent slower in Go
than C.  This compares very favorably to almost any other high level programming language.
So the author considered a modification to the algorithm above that could effectively use parallelism.

The goal was to modify the old method to a steepest descent search.
Instead of iterating through
the edges and immediately changing the color of edges whose recoloring reduced the number of
monochromatic subgraphs, the modified algorithm finds the edge that makes the
largest improvement and recolors that single edge.
This allows the $ecolor$ functions to be run in parallel.
Of course this change might be expected to make the algorithm less efficient, since it examines
every edge before making a single color change.  And this is indeed the case,
at least in terms of total CPU time used.  But the hope was that it would explore the search space
differently, and perhaps find different graphs than the old method found.
This was indeed the case, at least for $R(5,6)$.

The program was run using $414$ threads (on a computer capable of running $512$ concurrently). Since $K_{58}$
has $1653$ edges, each thread had at most four edges to consider.

\vspace{5mm}

\begin{algorithm}[H]
\KwData{$A$ is the $n \times n$ $2$-color adjacency matrix}
\KwData{$L$ is a list of all edges of $K_n$}
\KwData{$weight$ is a list of two elements $[5, 4]$ indexed by color}
\Begin{
  Randomly color $A$ using colors $1$ and $2$ with $Prob(color\,1) = 0.47$\;
  \Repeat{$ count(A)=0 $}  {
    \For{$e=uv$ in $L$}{
       $c \leftarrow a[u,v]$ \;
       $ \Delta(e) = weight[3-c] \cdot escore(e,3-c) - weight[c] \cdot escore(e,c) $\;
    }
    Change the color of the edge $e$ with minimum $\Delta(e)$ \;
    In case of ties, choose $e$ uniformly randomly \;
  }
}
\caption{A search procedure for finding $(5,6)$-colorings of $K_{n}$.  Note that if an edge has color $c$ then
$3-c$ is the other color.}
\end{algorithm}

\vspace{2mm}

Several different $(5,6)$-colorings of $K_{58}$ were found.  The one given below
was the first coloring found, and the only coloring found more than once.
An adjacency matrix of the new $(5,6)$-coloring is available at the
following location.

\begin{center}
\verb+http://cs.indstate.edu/ge/RAMSEY/g56.txt+
\end{center}

\noindent
The graph of the $K_5$-free color is given below in Brendan McKay's {\it graph6} format \cite{nauty}.

\begin{figure}
\begin{verbatim}
y?CbSgfEqHPAaVPjuDZSakemTaZl\tqzJctGuMgFWM`Z`[bEtB
~HyK~BYVtPNAyiPj}IuIMfwkQ}WnYIZfgtoVBiV[XVS]jUAlxa
qwRne_\n[al\Z]`Op~VNtgmlK|jsEY[ZyouKuUo\n]WVb`[vZD
ttdFpV^oHk[{pjjbfJLKfll[}mVHwIyVWkZ[sjw\qgxuYTqr\P
iWTa}xvQIqdDfjafWm~gT|xnJgN^KoZX{xOv~C@}|IizEJf^NY
sl[wUwnqVss~A[z}lY}FNqJrNv
\end{verbatim}
\caption{The graph6 representation of the color $1$ graph in a $(5,6)$-coloring of $K_{58}$.}
\end{figure}

This new variant of our search algorithm was also tried on other Ramsey problems,
particularly $R(5,5)$, $R(4,6)$, $R(3,10)$, and $R(3,3,3,3)$,
where it failed to find any colorings that improve the known lower bounds.  In each case, however,
it was able to equal the old bounds, but with known colorings.  One note of optimism:
it does seem to find the known colorings (of $K_{50}$) for $R(3,3,3,3)$ more easily than
other search methods this author has tried.

\bibliographystyle{plain}

\end{document}